\documentclass[12pt, twoside, leqno]{article}

\usepackage{amsmath,amsthm}
\usepackage{amssymb}
\usepackage{mathtools}
\usepackage[shortlabels]{enumitem}

\usepackage{enumitem}

\usepackage{graphicx}
\DeclareUnicodeCharacter{02B9}{'}

\usepackage[T1]{fontenc}

\usepackage{nowidow}

\usepackage{mathrsfs}

\newtheorem{theorem}{Theorem}[section]

\newtheorem{lemma}[theorem]{Lemma}

\theoremstyle{definition}

\newtheorem{remark}[theorem]{Remark}
\newtheorem{example}[theorem]{Example}

\numberwithin{equation}{section}

\newcommand{\ve}{\varepsilon}

\newcommand{\diam}{{\rm diam}}
\newcommand{\ra}{\rightarrow}

\newcommand{\bb}[1]{\mathbb{#1}}

\newcommand{\R}{\bb{R}}

\newcommand{\N}{\bb{N}}

\newcommand{\Oset}{\varnothing}

\usepackage[utf8]{inputenc}

\begin{document}


\baselineskip=17pt


\title{Revisiting the Lavrentiev Phenomenon  in~One Dimension}

\author{Wiktor Wichrowski\\
University of Warsaw \\ 
Faculty of Mathematics, Informatics
and Mechanics\\
Banacha 2, 02-097 Warszawa, Poland\\
E-mail: ww439108@students.mimuw.edu.pl
}

\date{}

\maketitle


\renewcommand{\thefootnote}{}



\renewcommand{\thefootnote}{\arabic{footnote}}
\setcounter{footnote}{0}


\begin{abstract}
We clarify and extend insights from Lavrentiev’s seminal paper. We examine the original theorem dealing with the absence of the Lavrentiev phenomenon, a cornerstone issue in the calculus of variations. We point out some inconsistencies in the original proof by providing a counterexample and supply the result with a new, concise, and complete reasoning. In the appendix, we also provide additional details to supplement the original proof. 
\end{abstract}

\section{Introduction}
Let us consider the variational functional
\begin{equation*}
	\mathcal{F}(u) = \int _a^b f\big(x, u(x), u'(x)\big)\,dx \,,
\end{equation*}
\noindent
where \(f\colon [a,b]\times \R \times \R \to \R \) and \(u\) belongs to some class \(\mathcal{E}\) (to be specified). The calculus of variations focuses on the study of such functionals, in particular on the problem of finding a function that minimizes \(\mathcal{F}\) within a given class of functions \(\mathcal{E}\), that is,
\[
\inf_{u \in \mathcal{E}} \mathcal{F}(u)\,.
\]
This involves investigating the properties of \(\mathcal{F}\) and of admissible functions \(u \in \mathcal{E}\), with the goal of determining whether a minimizer exists and, if so, characterizing such minimizers. Note that the infimum is not always attained within a certain class of functions. A thorough discussion of such situations, with illustrative examples, can be found in the work of Buttazzo, Hildebrandt, and Giaquinta \cite{buttazzo1998one}. In these instances, the infimum cannot be reached by restricting to regular functions alone, which highlights the necessity of considering broader classes of functions to obtain a reachable minimizer.

We would want to calculate the infimum by restricting to regular functions.
This raises an important question: under what conditions is such a restriction sufficient? Clearly, the following inequality always holds:
\[ \inf_{\text{all } u} \mathcal{F}(u) \leq \inf_{\text{regular } u} \mathcal{F}(u)\, .\]
However, there are situations in which there exists a constant \( c \in \mathbb{R} \) such that:
\[ \inf_{\text{all } u} \mathcal{F}(u) <c< \inf_{\text{regular } u} \mathcal{F}(u) \,.\]
Such a situation is called the Lavrentiev phenomenon after the prominent contribution by Michail Lavrentiev\footnote{
The most widely used Latin spelling, Lavrentiev, does not match the spelling used in the author’s original French publication, where his name appears as “M. Lavrientieff”. His name has appeared in a variety of other forms in the literature, including Lavrientieff, Lavrientief, Lavrentʹev, Lavrentjev, Lavrientiev, and Lavrentyev.
}
in \cite{Mascon}. Therein the first example of the occurrence of the phenomenon is provided together with the conditions necessary for its absence.  The phenomenon was later popularized by the example by M\'ania \cite{Mania}, the proof of which was later simplified (see e.g.~\cite{MR879878}). A new example of functional for which Lavrentiev phenomenon occurs has been presented by Cerf and Mariconda in \cite{cerf2024lavrentievphenomenon},
under the assumption that only one endpoint value of \(u\) is prescribed.
assuming only one boundary condition. Other examples can be found in \cite{Bebugen}. In context of the one-dimensional Lavrentiev phenomenon, a topic was thoroughly examined in the book by Buttazzo, Hildebrandt, and Giaquinta \cite{buttazzo1998one}. Mariconda \cite{carlo} explores its emergence in Sobolev spaces, identifying conditions under which it can be avoided, particularly through boundary considerations, while Gómez \cite{gomez} extends the analysis to capillarity-type problems and solution regularity in bounded variation spaces. Furthermore, Buttazzo and Mizel \cite{VictorMizel} showed that the Lavrentiev phenomenon is local in nature.

The exclusion of the Lavrentiev phenomenon is significant in the regularity theory of minimizers to variational functionals. 
In the case of non-autonomous functionals with non-standard growth conditions, the presence of the Lavrentiev phenomenon can obstruct the regularity of minimizers and their gradients.
Recent work has made advances in this area, with studies on double-phase integrals, \((p, q)\)-growth conditions, and higher integrability of minima, as well as comprehensive reviews of these complex variational problems \cite{BalciLars, Baroni, Filippis, Eleuteri,  Esposito,Mingione,Zhikov, Zhikov2}. Notably, absence of the Lavrentiev phenomenon’s was established without relying on growth conditions or convexity in \cite{Bousquet2023, Bousquet}. New studies identify conditions preventing the Lavrentiev gap in both scalar and vectorial cases in the anisotropic setting, under certain balance conditions \cite{Borowski3, Borowski2, Borowski1, koch2023lavrentievgapconvexvectorial}. Results concerning the absence of the Lavrentiev phenomenon are crucial in researching the stability and convergence of numerical methods. This phenomenon can lead to failures in traditional numerical techniques, such as finite element methods, which may struggle to detect singular minimizers. This issue was studied in the one-dimensional setting by Ball and Knowles~\cite{Ball_Knowles}, which examined Manià's example to illustrate the phenomenon's impact on numerical methods. Similar challenges in higher dimensions are addressed in \cite{Bai_Zhi, Li_Zhi, Sivaloganathan, XU-2011}, which propose more sophisticated numerical methods for simulating the rapid expansion of voids and other complex deformations under tensile stress.

In this paper, we revisit Lavrentiev’s seminal 1927 work \cite{Mascon}. Specifically, we examine his original example, including a review of the flaws in the original proof and the presentation of a new proof  (Example~\ref{ex:lav_org}). Following this, we focus on the proof of the absence of the Lavrentiev phenomenon as a consequence of his Approximation Lemma. We review Lavrentiev’s original proof from his seminal paper, refining the details and presenting a simple counterexample to demonstrate that certain aspects are not immediately obvious. Subsequently, in Section~\ref{sec:proof}, we present a new concise, and complete proof of Lavrentiev’s Approximation Lemma. In Appendix~\ref{Appendix_A} we correct  Lavrentiev's original proof. Specifically, we explain how to construct the sets required for his argument, as the details provided in his paper are misleading. 

\bigskip

\noindent \textbf{Notation. }  We define \(AC_*[a, b]\) as the set of absolutely continuous functions on the interval \([a, b]\) that satisfy the boundary conditions \(x(a) = A\) and \(x(b) = B\):
    \[
        AC_*[a, b] = \{ x \in AC([a, b]) \colon x(a) = A, x(b) = B \}\,.
    \]
    We define \(Lip_*([a, b])\) and \(C^k_*([a, b])\) for \(k = 1, \ldots, \infty\), to denote functions in the respective spaces that also satisfy these boundary conditions.

      A set is called \emph{perfect} if it is closed and contains no isolated points. 
\section{The Lavrentiev Example}
In his original paper, namely \cite{Mascon}, Michail Lavrentiev presented an implicit example of a functional, for which the Lavrentiev phenomenon occurs, and a complicated justification of this fact. He postulated the existence of a smooth function \( f\colon \mathbb{R} \to \mathbb{R} \) satisfying certain assumptions, for which a gap occurs. His functional took the following form:
\begin{equation}\label{eq:the-funtional}
    \mathcal{F}(u) \coloneqq \int_0^1 e^{-\frac{2}{(u(x) - \sqrt{x})^2}} f(u'(x)) \, dx, \quad \text{with} \quad u(0) = 0  \text{ and } u(1) = 1\,.
\end{equation}
We will outline the original proof from  \cite{Mascon}, show some of its flaws, and present new elementary reasoning instead.

\begin{figure}[htbp]
    \centering
    \includegraphics[width=0.732\linewidth]{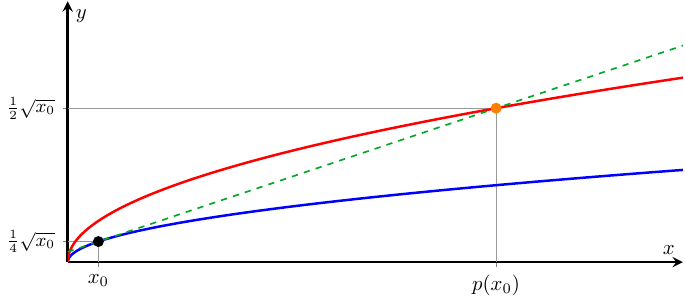}
    \caption{\small \textit{
        \begin{tabular}{@{} l l @{}}
            In blue: & Plot of $y=\frac{1}{4}\sqrt{x}$.\\
            In red: & Plot of $y=\frac{1}{2}\sqrt{x}$.\\
            In green: & Tangent line to plot of $y=\frac{1}{4}\sqrt{x}$ in point \(x_0\).
        \end{tabular}}}
    \label{fig:ex_lav_p}
\end{figure}

Lavrentiev required the function \( f \) to have a particular set of properties. To introduce these conditions, we first define the auxiliary function \( p\colon [0,\infty) \to \mathbb{R} \) as follows
\begin{equation}\label{lav:def_p}
    p(x) \coloneqq 7x + 4\sqrt{3} x\,.
\end{equation}
This function arises from the following geometric construction. For each \( x_0 \in \mathbb{R} \), we consider the tangent line to the parabola \( y(x) = \frac{1}{4} \sqrt{x} \) at the point \( \left( x_0, \frac{1}{4} \sqrt{x_0} \right) \). The value \( p(x_0) \) represents the abscissa (i.e., the \( x \)-coordinate) of the second intersection point of this tangent line with the parabola \( y = \frac{1}{2} \sqrt{x} \), moving from left to right (see Figure~\ref{fig:ex_lav_p}).
Now we are ready to present the properties outlined by Lavrentiev in his paper. Some of them  have been reformulated for improved readability. 
\noindent
First, Lavrentiev required that \( f \) is increasing and convex with a minimum at \( 0 \) that is greater than \( 1 \):\begin{equation}\tag{I}\label{I}
    \frac{d^2 f}{dx^2} > 0 \quad\text{ and } \quad \min_{x \in \mathbb{R}} f(x) = f(0)\geq 1\,.
\end{equation}
\noindent 
The second property yields a lower bound for the growth of \( f \), outside of a neighborhood of~\( 0 \):
\begin{equation}\tag{II}\label{II}
\inf_{\substack{\frac{1}{2} \leq a < 1 \\ p_0 \in (a, 1]}}
\left( {(p_0 - a)} f \left[ \frac{ \frac{1}{4} \sqrt{\frac{1}{2}} }{p_0 - a} \right] \right) > e^{8\sqrt{2}}\,.
\end{equation}
\noindent
Additionally, Lavrentiev specified that the growth of $f$ should satisfy the following condition across the entire domain:
\begin{equation}\tag{III}\label{III}
e^{-\frac{8}{\sqrt{x}}} \big(p(x) - x\big) f \left[ \frac{1}{4} \left( \frac{\sqrt{p(x)} - \sqrt{x}}{p(x) - x} \right) \right] > 1 \quad\quad (0 < x \leq 1)\,.
\end{equation}
%
\noindent
The fourth condition ensures that \( f \) is bounded from below:
\begin{equation}\tag{IV}\label{IV}
\frac{e^{-\frac{8}{\sqrt{x}}}}{4} \sqrt{x} f\left( \frac{1}{8 \sqrt{x}} \right)
\cdot \frac{1}{\sqrt{1 + \left(\frac{1}{8 \sqrt{x}}\right)^2}} > 1 \quad \quad (0 < x < 1)\,.
\end{equation}
\noindent Finally, Lavrentiev's additional monotonicity condition: \begin{equation}\tag{V}\label{V}\frac{d}{dt} \left( f(t)\cdot \frac{1}{\sqrt{1 + t^2}} \right) > 0 \quad \quad  (t > 0)\,.\end{equation}
He did not provide complete proof for the existence of a function satisfying the above properties, offering only brief remarks. Before presenting the original reasoning, we will need the following lemma, which follows from Jensen's inequality.
\begin{lemma}\label{classic}
    The minimum of 
    \[
    \mathcal{F}(u) \coloneqq \int_a^b f(u'(x)) \, dx\,,
    \]
    for \(u \in \text{Lip}_*([a,b])\) and any convex function \(f\colon \mathbb{R} \to \mathbb{R}\) is attained when \(u(x)\) is a straight line. Therefore, the following inequality holds:
    \[
    \mathcal{F}(u) \geq (b-a) f\left( \frac{u(b) - u(a)}{b-a} \right).
    \]
\end{lemma}
\noindent Let us now examine the original proof, which involves assumptions \eqref{I}–\eqref{V}.
\newpage
\textbf{Sketch of Lavrentiev's Original Proof}
\begin{enumerate}[Step 1.]
    \item Define \(\varphi\) as the infimal value of functional $\mathcal{F}$ from \eqref{eq:the-funtional} on the interval from $a$ -- the abscissa of the first point of intersection of the graph of $u$ with the graph of \( y(x) = \frac{\sqrt{x}}{4} \) to $1$. Namely,\footnote{The notation has been modified from the original: \( a \coloneqq x_0 \), \( b \coloneqq x' \), \( u(x) \coloneqq y(x) \), and \(\varphi(a) \coloneqq \varphi(a, \frac{1}{4} \sqrt{a})\).}
    \begin{equation*}
        \begin{aligned}
            \varphi(a) \coloneqq \inf \Bigg\{ \int_{a}^{1}&e^{-\frac{2}{(u(x) - \sqrt{x})^2}} f(u'(x)) \, dx \colon
            \\ 
            & u \in C^1([a, 1]), \, u(a) = \frac{1}{4} \sqrt{a}, \, u(1) = 1 \Bigg\}\,.
        \end{aligned}
    \end{equation*}
    Note that the integral above is taken over the interval \( [a, 1] \) instead of \( [0, 1] \). 
    Consequently, it is sufficient to prove that \(\varphi(a) > 1\) for all \( a \in (0,1) \).
    
    \item Define a sequence \(\{x_n\}_{n=1}^\infty\) recursively by setting \(x_1 = \frac{1}{2}\), and for \(n \geq 1\), let \(x_{n+1}\) satisfy the equation \(p(x_{n+1}) = x_n\), where the function \(p\) is defined as in \eqref{lav:def_p}. Notice that \(x_n \to 0\) as \(n \to \infty\). Therefore, it is sufficient to prove that \(\varphi(a) > 1\) for all \(a \geq x_n\) for all \(n \in \mathbb{N}\). This can be accomplished by induction.
    \item \textbf{Base Case:} Conditions \eqref{I} and \eqref{II}, combined with Lemma~\ref{classic}, imply that \(\varphi(a) > 1 \) for all \(a > x_1=\frac{1}{2}\). 
    \item \textbf{Induction Step:} Assume that the inequality holds for all \(a \geq x_n\), and show that \(\varphi(a) > 1\) holds for \(a \geq x_{n+1}\).  To demonstrate this, consider the following three cases (see Figure~\ref{fig:lav_fig}), where $y_1=u(p(a))$:
        \begin{enumerate}[{Case} 1)]
            \item \(y_1 \leq \frac{1}{4} \sqrt{p(a)}\). This case is straightforward because \(x_n < p(a)\), allowing us to use our induction assumption. Indeed, in such a case there always will be a point \(\bar{a} \geq x_n\) at which the graphs of  \(u\) and \(x\mapsto\frac{1}{4}\sqrt{x}\)  intersect.
            \item \(\frac{1}{4} \sqrt{p(a)} \leq y_1 \leq \frac{1}{2} \sqrt{p(a)}\), and \(u(\xi) < \frac{1}{2} \sqrt{\xi}\) for \(a \leq \xi \leq p(a)\). This is the case where \(u(x)\) is separated from \(\sqrt{x}\) over a long interval. It is slightly more involved than the first case, as it requires the use of conditions \eqref{I} and \eqref{III}.
            \item \( y_1 > \frac{1}{2} \sqrt{b} \), for \( a < b < p(a) \), meaning the curve \( u = u(x) \) intersects the parabola \( u = \frac{1}{2} \sqrt{x} \) at least once. This case is the most complicated and contains some flaws. The author attempted to use  condition \eqref{IV} and \eqref{V}. He claimed that:
            \[
            \frac{\frac{1}{2}\sqrt{b} - \frac{1}{4}\sqrt{a}}{b - a} \geq \frac{1}{8\sqrt{a}}\,,
            \]
            which can be written as:
            \begin{equation*} 
            14ab \geq b^2 + a^2\,.
            \end{equation*}
            Notice that this cannot hold in general.     As an example, when \( b = 10a \), we have \( 10a < 13.92a \approx p(a) \) and the above fails to hold.
        \end{enumerate}
\begin{figure}[htbp]
    \centering
    \includegraphics[width=0.9\linewidth]{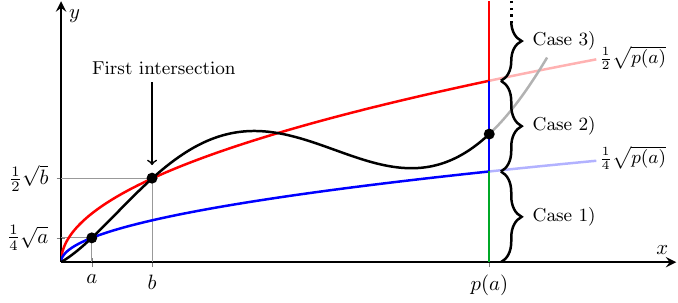}
    \caption{\small \textit{
        \begin{tabular}{@{}l@{\hskip 18pt}l@{}}
            In blue: & Plot of \(y=\frac{1}{4}\sqrt{x}\).\\
            In red: & Plot of \(y=\frac{1}{2}\sqrt{x}\).\\
            In green: & 
            Example of \(y(\cdot)\), with the first intersection point marked.
            \vspace{4pt}\\
            & {\hspace{-2.17cm}The vertical line at \(p(x_0)\) indicates cases: for functions that do not  
            } \\
            & {\hspace{-2.17cm}intersect with \(y=\frac{1}{2}\sqrt{x}\), we focus on the intersection point relative to  
            } \\
            & {\hspace{-2.17cm}the position of \(y=\frac{1}{4}\sqrt{p(x_0)}\).
            } \\
        \end{tabular}
    }}
    \label{fig:lav_fig}
\end{figure}
\end{enumerate}
We now present a new, alternative proof for this example.
\begin{example}[Lavrentiev 1926]\label{ex:lav_org}
	Consider the functional \(\mathcal{F}\) defined by
	\begin{equation}\label{Lav}
		\mathcal{F}(u) \coloneqq \int_0^1 e^{-\frac{2}{(u(x) - \sqrt{x})^2} } f(u'(x)) \, dx, \quad \text{with} \quad u(0) = 0 \quad \text{and} \quad u(1) = 1\,.
	\end{equation}
	Then, there exists a smooth function \(f\colon \R \to \R \) such that
 the Lavrentiev phenomenon occurs between \(AC_*([0,1])\) and \(C^1_*([0,1])\), namely:
	\[
	\inf_{u \in AC_*([0,1])} \mathcal{F}(u) = 0 \quad \text{and} \quad \inf_{u \in C^1_*([0,1])} \mathcal{F}(u) > 0\,.
	\]
\end{example}
\begin{proof}  
Consider the sequence of functions \(u_n(x) \coloneqq \sqrt{x} + \frac{x(x-1)}{n}\) for $n\in\N$. For this sequence, we have{

	\begin{align*}
	\inf_{u \in AC_*([0,1])} \mathcal{F}(u) = \lim_{n\to \infty}\mathcal{F}(u_n)&\leq \lim_{n\to \infty } \int_{\tfrac{1}{4}}^{\tfrac{3}{4}}\left|  e^{-\frac{2n^2}{x^2(1-x)^2}}f(u_n'(x))\right| \, dx \\ &\leq \lim_{n\to \infty } C \int_{\tfrac{1}{4}}^{\tfrac{3}{4}} e^{-n^2}\, dx = 0\,,
		\end{align*}
}
indicating that the infimum over absolutely continuous functions is indeed \(0\). {
We demonstrate that, for the functional $\mathcal{F}$ with the choice of
    \begin{equation}\label{lav:def_f}
        f(x)\coloneqq  e^{8\cdot256x^2}\,,
    \end{equation}
    we have
    \[
	\inf_{u \in C^1_*([0,1])} \mathcal{F}(u) > e\,.
	\]
}

Assume \( u \in \text{Lip}_*([0,1]) \). Its derivative \( u'(x) \) is bounded by a Lipschitz constant \( L \) for all \( x \in [0,1] \). Since the derivative of the function \( x \mapsto \sqrt{x} \) becomes unbounded as \( x \) approaches \( 0 \), we can conclude that there exist \( a, b \in (0,1) \) such that
\begin{equation}\label{def:ux}
    \frac{\sqrt{x}}{4} < u(x) < \frac{\sqrt{x}}{2} \quad \text{for all } x \in (a, b)\,,
\end{equation}
and equality on boundary, i.e.:
\[u(a) = \frac{\sqrt{a}}{4} \quad \text{ and } \quad u(b) = \frac{\sqrt{b}}{2}\,. \]
\noindent Thus, for all \( x \in [a, b] \), the following inequality holds:
\begin{equation*}
    e^{-\frac{2}{(u(x) - \sqrt{x})^2}} \geq e^{-\frac{2}{\left(\frac{\sqrt{x}}{2} - \sqrt{x}\right)^2}} = e^{-\frac{8}{x}}\,.
\end{equation*}

\noindent We will divide this proof into two cases.\\
\noindent
\textbf{Case 1: } \(b<2a\). 
Young's inequality can be written as
\[
 \beta^2\ve \geq 2\alpha\beta - \frac{\alpha^2}{\ve}\,,
\]
so by substituting \( \ve \coloneqq e^{-\frac{8}{x}} \), \( \alpha \coloneqq 1 \), and \( \beta^2 \coloneqq f(u'(x)) \), we obtain
\[
f(u'(x)) e^{-\frac{8}{x}} \geq 2 \sqrt{f(u'(x))} - e^{\frac{8}{x}} \,.
\]
\noindent
Integrating this inequality over the interval \([a,b]\), we get
\begin{equation}\label{lav:int}
    \int_a^b f(u'(x))e^{-\frac{8}{x}} \, dx \geq \int_a^b \left( 2 \sqrt{f(u'(x))} - e^{\frac{8}{x}} \right) \, dx \,.
\end{equation}
\noindent
From the assumption of Case 1, that is \( 0 < a < b < 2a \), and since \(\sqrt{2}+\frac{1}{2} <2\), we have the following inequality
\begin{equation}\label{lav:class}
\begin{aligned}
        \int_{a}^{b}\left( u'(x)\right)
        = &\, (b-a) \left( \frac{\frac{1}{2}\sqrt{b}-\frac{1}{4}\sqrt{a}}{b-a}\right)\\
        = & \, (b-a) \left( \frac{b-\frac{1}{4}a}{2(b-a)( \sqrt{b}+\frac{1}{2}\sqrt{a})}\right)\\
        \geq &\,  (b-a) \left( \frac{\frac{3}{16}\sqrt{a}}{b-a}\right) = \frac{3}{16}\sqrt{a}\,.
\end{aligned}
\end{equation}
\noindent For convenience, we define
\begin{equation}\label{lav:def_g}
    g(x) \coloneqq \sqrt{f(x)} = e^{4 \cdot 256 x^2}\,,
\end{equation}
which is an increasing and convex function. We now apply  \eqref{lav:int}  along with Jensen's inequality and the estimates from \eqref{lav:class}. This yields the following chain of inequalities:
\begin{equation*}
\begin{aligned}
    \int_a^b f(u'(x)) e^{-\frac{8}{x}} \, dx 
    & \stackrel{\eqref{lav:int}}{\geq} \int_a^b \left( 2\sqrt{f(u'(x))} - e^{\frac{8}{x}} \right) \, dx \\
    &\;\; \geq (b - a) \left( \frac{1}{b - a} \int_a^b 2\sqrt{f(u'(x))} \, dx - e^{\frac{8}{a}} \right) \\
    & \stackrel{\eqref{lav:def_g}}{=} (b - a) \left( \frac{1}{b - a} \int_a^b 2g\left(u'(x)\right) \, dx - e^{\frac{8}{a}} \right) \\
    & \stackrel{\text{Jensen}}{\geq} (b - a) \left( 2g\left( \frac{1}{b - a} \int_a^b u'(x)\, dx \right) - e^{\frac{8}{a}} \right) \\
    & \stackrel{\eqref{lav:class}}{\geq} (b - a) \left( 2g\left(
     \frac{\left(\frac{3}{16}\right)\sqrt{a}}{b-a}
    \right) - e^{\frac{8}{a}} \right).
\end{aligned}
\end{equation*}
By the definition of \(g\) given in \eqref{lav:def_g}, we can estimate the integral \eqref{Lav} under the assumptions of Case 1, i.e., \(0 < a < b < 2a\):
\begin{equation*}
    \begin{aligned}
        \int_a^b f(u'(x)) e^{-\frac{8}{x}} \, dx 
        & \geq (b - a) \left( 2e^{ \frac{8a}{(b - a)^2}} - e^{\frac{8}{a}} \right) \\
        & \geq (b - a) \left( e^{ \frac{8a}{(b - a)^2}} + e^{\frac{8}{a}} - e^{\frac{8}{a}} \right) \\
        & \geq (b - a) e^{ \frac{1}{b - a}} \geq e\,.
    \end{aligned}
\end{equation*}
Last inequality follows from fact that \(x\mapsto xe^{\frac{1}{x}}\) is decreasing for \(x \in(0,1]\). This concludes the proof in the first case.\\
\noindent
\textbf{Case 2:} \(b \geq 2a\). {Keeping in mind \eqref{def:ux} we estimate
\begin{equation*}
\begin{aligned}
    \int_a^b f(u'(x)) e^{-\frac{8}{x}} \, dx 
     \geq& e^{-\frac{8}{a}} \int_a^{2a} f(u'(x)) \, dx
     \\
    \stackrel{\text{Lemma }\ref{classic}}{\geq} (2a&-a) e^{-\frac{8}{a}} f\left(\frac{u(2a)- \frac{1}{4}\sqrt{a}}{2a - a}\right) 
    \\
     \geq  ae^{-\frac{8}{a}}& f\left(\frac{\sqrt{2}-1}{4\sqrt{a}}\right) 
     =   a e^{\frac{8\cdot16 (3-2\sqrt{2})-8}{a }}\geq 
     e
\end{aligned}
\end{equation*}
The last inequality holds, since \(a \geq e^{ -\frac{1}{a}}\), for \(a>0\). This concludes the proof in the second case. \newline}

By combining both cases, we establish that the Lavrentiev phenomenon occurs between \(AC_*([0,1])\) and \(C^1_*([0,1])\) for the functional \eqref{Lav}, with \(f\) defined as in \eqref{lav:def_f}.
\end{proof}
\section{Absence of the Lavrentiev Phenomenon in one-dimensional problems}

In 1927 Lavrentiev published an article in which he formulated conditions sufficient for the absence of a gap between absolutely continuous functions and \(C^1\). Firstly, we will present an approximation lemma and theorems resulting from it. Then, we will discuss the original proof of the lemma and provide a simplified proof.

\begin{lemma}[Approximation Lemma, \cite{Mascon}] \label{Lav1}
 Let  {$f\colon [0,1]\times\R\times \R\to[0,\infty)$   be a continuous function of variables $(x,y,\xi )$} and such that $\left|\frac{\partial f}{\partial y}\right|< M$. Let $u \in AC_*([0,1])$ such that the function $x\mapsto f[x, u(x), u'(x)]$ is integrable. Then, for  {every} $ \varepsilon > 0$, there exists a function  \(\phi \in C^\infty ([0,1])\) satisfying the following conditions:
		\begin{enumerate}[(L1)]
			\item\label{Lav1:1}  $\phi(0) = u(0)$ and $\phi(1) = u(1)\,,$
			\item\label{Lav1:2}  $|u(x) - \phi(x)| < \varepsilon \,$ for all \(x\in [0,1]\,,\)
			\item\label{Lav1:3}  \( \left| \int _0 ^1 f(x,u(x),u'(x))\, dx -\int _0 ^1 f(x,\phi(x),\phi'(x))\, dx \right| < \varepsilon\,. \)
		\end{enumerate}
\end{lemma}

Thanks to Lemma~\ref{Lav1}, we can exclude Lavrentiev's phenomenon for functions that have bounded derivatives with respect to the second variable. Thus, we have
\[
\inf_{AC_*([0,1])}\mathcal{F} = \inf_{C^1_*([0,1])}\mathcal{F}\,.
\]
This is the content of Theorem~\ref{coro2}.

In his 1927 paper, to prove Lemma~\ref{Lav1}, Lavrentiev constructed a piece-wise linear function that satisfied the conditions~\ref{Lav1:1}-\ref{Lav1:3} from the lemma.  While the original formulation of Lavrentiev's Approximation Lemma is true, there is a step in his proof that is not adequately demonstrated. Let us examine this issue.

\subsection{Sketch of original proof}

The original proof of Lavrentiev's Approximation Lemma from \cite{Mascon} revolves around constructing a piece-wise constant function that approximates \(u\) in a suitable manner. Here's a sketch of the original proof:

\begin{enumerate}[1.]
	\item Notice that it suffices to prove that there exists a piece-wise linear function satisfying the same conditions as the function \(\phi\).
	      
	\item Let \(\ve>0\) be given and use Lusin's Theorem to take take \(P_\varepsilon\) such that:
	      \begin{enumerate}
	      	\item \(u'\) is continuous on it, 
	      	\item the measure of \(P_\varepsilon\) approaches 1 as \(\varepsilon\) approaches 0, 
	      	\item \(P_\varepsilon\) is perfect -- it is closed and contains no isolated points.
	      \end{enumerate}
	      
	\item  \label{org:step3} Take a finite partition \(\mathcal{A} = \{A_1, A_2,\ldots, A_n \}\) of \(P_\varepsilon\) such that
	      \[
	      	\forall_{i\in \{1,2,\hdots, n\}}\quad \exists_{c_i\in \R} \quad \forall_{x\in A_i} \quad |u'(x) - c_i| < \varepsilon\,.
	      \]
	\item \label{org:step4} Expand \(\mathcal{A}\) to the family \(\mathcal{B} = \{B_1, B_2, \ldots, B_n\}\) of disjoint intervals such that for all \(i \in \{1, 2, \ldots, n\}\), we have \(A_i \subseteq B_i\). Take the complement \(B_i \setminus A_i\), divide it into intervals, and denote them by family \(\left\{C_{i,k}\right\}_{k=1}^{\infty}\).
  
	\item Construct a piece-wise constant function
	      \begin{equation*}
	      	\psi(x) \coloneqq
	      	\begin{cases}
	      		c_i & \text{if } x \in A_i \text{ or } x \in C_{i,k} \text{ for } k > k_i \,,\\
	      		0   & \text{otherwise,}                                                    
	      	\end{cases}
	      \end{equation*}
	      where \(\sum\limits _{k=k_i}^\infty \lambda \left( C_{i,k}\right) \) is sufficiently small.
	\item Construct piece-wise linear function
	      \begin{equation*}
	      	\bar{u}(x) \coloneqq \int _0 ^x \psi(t)\,dt+u(0)
	      \end{equation*}
	      and show that it satisfies conditions~\ref{Lav1:2} and~\ref{Lav1:3} from Lemma~\ref{Lav1}.

	\item Fix the boundary condition: to do this, choose a subset \(E\) of \(P_\varepsilon\) of measure \(1/2\)   and construct the function accordingly, similarly as in Steps 5 and 6{, set} 
 
	      \begin{equation*}
	      	\bar{\psi}(x) \coloneqq
	      	\begin{cases}
	      		\psi(x) +2\left( \bar{u}\left(1\right) - u\left(1\right) \right) & \text{for } x \in E \,,\\
	      		\psi(x)                                   & \text{otherwise,}    
	      	\end{cases}
	      \end{equation*}
       {and }finally, define
	      \begin{equation*}
	      	\bar{\bar{u}}(x) \coloneqq \int_0^x \bar{\psi}(t) \, dt + u(0)\,.
	      \end{equation*}
       
	\noindent
    Verify that this function satisfies conditions~\ref{Lav1:1}--\ref{Lav1:3}.
\end{enumerate}

\subsection{Simple counterexample}\label{Simple_counterexample}
\noindent In the original proof, it is not clarified why we can choose a family from step \ref{org:step4} Firstly, let us consider a simple example showing it is not obvious that those sets exist.

\begin{lemma}
    For every \(\ve \in (0,\tfrac{1}{4})\), there exists a set \(P_\ve\) and \(u\in AC_*([0,1])\) which has the following properties:
    \begin{enumerate}
        \item \(P_\ve\) is perfect (closed and without isolated points),
        \item \(\lambda(P_\ve) > 1 - \ve\),
        \item \(u'\) is continuous on \(P_\ve\),
        \item it is not possible to partition \(P_\ve\) into a finite number of intervals.
    \end{enumerate}
\end{lemma}
\begin{proof} The idea is to take \(P_\ve\) that contains infinitely many connected components.  Take any \(\ve \in (0,\tfrac{1}{4})\,.\) We define \(u'(x)\) as follows

\begin{equation*}
    u'(x) \coloneqq \begin{cases}
        -1 & \text{ if } x \in \left[\frac{1}{2^{2n+1}}, \frac{1}{2^{2n}}\right],\ n \in \mathbb{N}\cup \{0\}\,,\\
        1 & \text{ otherwise.}\\
    \end{cases}
\end{equation*}
Consider the set  \(P_\ve \coloneqq \bigcup_{n=1}^\infty \left[\frac{1+\ve}{2^{n}}, \frac{1-\ve}{2^{n-1}} \right] \cup \{0\}\). This set is perfect, and its measure is \(1-3\ve\). Additionally, \(u'\) is continuous on \(P_\ve\).

\begin{figure}[htbp!]
    \centering
     \includegraphics[width=0.9\linewidth]{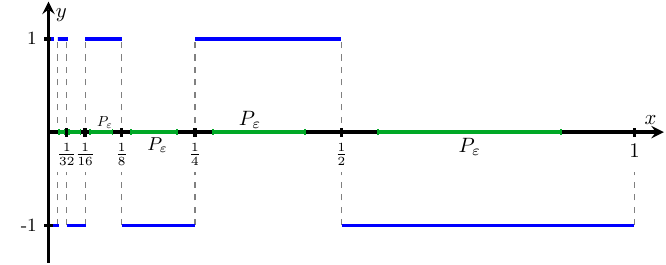}
    \caption{In blue, there is graph of \(u'\), while on the x-axis in green the set \(P_\varepsilon\) is marked.}
\end{figure}
 
For any partition \(\mathcal{B} = \{ B_1, B_2, \ldots, B_N\}\) of \([0,1]\) into \(N\in \mathbb{N}\) intervals, there exists a set \(B_i\) containing \(x,y\in P_\ve\) such that \(u'(x)=1\) and \(u'(y)=-1\). That  contradicts the assumption that we can expand the division of \(P_\ve\) into a finite number of intervals piece-wise disjoint.
\end{proof}

In this example, \(P_{\ve}\) is specifically chosen to demonstrate that Lavrentiev's argument is not correct. This raises the question of whether there exists a function for which there is no partition that can be extended to a finite number of closed, non-intersecting intervals.
This problem is discussed in detail in Appendix~\ref{Appendix_A}, where we present method for construction of those intervals.
\section{A Simplified Proof of the Lavrentiev Approximation Lemma}\label{sec:proof}

Here, we present a new, concise proof of Lemma~\ref{Lav1}. The idea of our proof is as follows. Since our function $u$ is assumed to be absolutely continuous (hence bounded), we can further assume without loss of generality that $0 \le u \le 1$. We will construct an approximate sequence in three steps: first, we approximate the function by Lipschitz functions; next, we smooth the approximation using the standard mollifier. Finally, we correct the end-point values. At each step, we ensure that the new function is sufficiently close to the previous one. Ultimately, we verify that the final function closely approximates \(u\) using the triangle inequality.

\begin{enumerate}[Step 1.]
	\item \label{step1}\textbf{Approximation by a sequence of Lipschitz functions (without keeping end-point values)}.\\
 Let us begin the construction from the derivative
      \begin{equation}\label{lav:v_k}
        v_k(x) \coloneqq
        \begin{cases}
            u'(x) & \text{if } |u'(x)| < k\,, \\
            0     & \text{otherwise.}       
        \end{cases}
      \end{equation}
      It is clear that $\{v_k\}_{n=1}^\infty$ converges a.e. to $u'$. Now, let
      \begin{equation*}
        u_k(x) \coloneqq \int_0^x v_k(t) \,dt \,.
      \end{equation*}
      Next, let's check that \(u_k \to u\) uniformly on \([0,1]\), that is
      \begin{equation}\label{step1:2}
        \forall _{\varepsilon>0} \quad \exists _k \quad \forall _x \quad \left|u_k(x) - u(x)\right| < \varepsilon \,.
      \end{equation}
      From the definition of \(AC([0,1])\) function, we get 
      \begin{equation*}\begin{aligned}
          |u_k(x) - u(x)| &= 
        \left| \int_{[0,x]\cap\{|u'(t)|>k\}}u'(t) \,dt \right|\\&\leq
        \int_{\{|u'(t)|>k\}} |u'(t)| \,dt \xrightarrow{k\ra \infty}0 \,.
      \end{aligned}
      \end{equation*}    
    Thus, condition~\ref{Lav1:2} is satisfied.
Now, let us check the third condition~\ref{Lav1:3}. We can infer that
\begin{equation}\label{step1:3}
    \left| \int_0^1 f(x,u(x),u'(x)) \,dx - \int_0^1 f(x,u_k(x),u_k'(x)) \,dx \right| \leq \varepsilon
\end{equation}
using the Dominated Convergence Theorem. Given  \eqref{lav:v_k} from definition of $u_k$ we obtain that \(u_k'(x) = 0\) or \(u_k'(x) = u'(x)\), so we have:
\begin{equation*}
\begin{aligned}
        \left| f(x,u_k(x),u_k'(x)) \right| \leq \big| f(x,u_k(x),u'(x)) \big| + \big| f(x& ,u_k(x),0) \big| \\ 
        &\eqqcolon G_k(x) + H_k(x)\,.
\end{aligned}
\end{equation*}

The first term can be written as
\[
G_k(x) \leq \left| f(x,u(x),u'(x)) \right| + M \big| u(x) - u_k(x) \big|\,,
\]
because we assumed that \( \left| \frac{\partial f}{\partial y} \right| < M \). This term is integrable since \( f(x,u(x), u'(x)) \) is integrable by our assumptions, and \( u(x) - u_k(x) \) is small due to the uniform convergence.

For the term \(H_k(x)\), notice that the sequence \( \{u_k\}_{k=1}^\infty \) is uniformly bounded (since it is uniformly convergent sequence of bounded real-valued functions), so we can approximate it by the supremum of \( f \) on a compact set. Thus, the assumptions of the Dominated Convergence Theorem are satisfied, and therefore \eqref{step1:3} is arbitrarily small for \( k \) large enough.

\item\label{step2}  \textbf{Approximation by smooth functions (without keeping end-point values).}\\
Let us apply mollification to the function \( u_k \). Define \(\eta_\varepsilon(x)\) as  
\[
\eta_\varepsilon(x) = \frac{1}{\varepsilon} \eta\left(\frac{x}{\varepsilon}\right),
\]  
where \(\eta\in C^\infty_c(\R)\) is the standard mollifier. Then, the mollification of \( u_k \) is defined as  
\[
u_{k,n} \coloneqq u_k \ast \eta_{\frac{1}{n}},
\]  
where \( \ast \) denotes the convolution operation.  { To define this properly, we extend \( u_k \) to a neighborhood of \([0,1]\) -- on the left by setting \( u_k(x) = u_k(0) \), and on the right by setting \( u_k(x) = u_k(1) \).} It follows that \( u_{k,n} \in C^{\infty}(\mathbb{R}) \) and that the sequence \( \{u_{k,n}\} \) uniformly converges to \( u_k \) as \( n \to \infty \) (see Appendix C.4 in \cite{evans10}).  Additionally, since \( u_k \) is Lipschitz with constant~\( k \) (see~\eqref{lav:v_k}), we obtain
    \begin{equation*}
        \begin{aligned}
            |u_{k,n}(x) - u_k(x)| 
                                    \leq&\int_{-1}^1 \left| (u_k(x - \tfrac{t}{n}) - u_k(x)) \right| \eta(t) \, dt     \\
                                    \leq&\int_{-1}^1 \left| k\,\frac{t}{n} \right| \eta(t) \, dt \leq \frac{k}{n}\,.
        \end{aligned}
    \end{equation*}
We will be checking conditions~\ref{Lav1:2} and~\ref{Lav1:3} between \(u_{k,n}\) and \(u_k\). In the last step, we will check the closeness of the approximation and the energy between \(u\) and the constructed function. We can see that condition~\ref{Lav1:2} is satisfied
    \begin{equation}\label{step2:2}
        \forall_{\varepsilon>0} \quad  \forall_{k }\quad \exists_{n_k} \quad\forall_{n\geq n_k} \quad \forall_{x} \quad |u_{k,n}(x) - u_k(x)| \leq \frac{k}{n} <\ve\,.
    \end{equation}

    Condition~\ref{Lav1:3} is slightly more involved; we need to use the Dominated Convergence Theorem. Firstly, observe that \(u_{k,n}\) are Lipschitz function with constant \(k\)
    \begin{multline*}
           |u_{k,n}(x_1)-u_{k,n}(x_2)|= \left|\int_{-1}^1 \big(u_{k}(x_1-y)-u_{k}(x_2-y) \big)\eta_{\frac{1}{n}}(y)\, dy \right| \\
           \leq \int_{-1}^1 \big|u_{k}(x_1-y)-u_{k}(x_2-y) \big|\eta_{\frac{1}{n}}(y)\, dy \leq k |x_1-x_2|\,.
    \end{multline*}

    Therefore a suitable bound is given by
    \begin{equation*}
        \begin{aligned}
              & \left| f(x, u_{k,n}(x), u'_{k,n}(x)) \right| \leq \sup_{(x,y,z)\in [0,1]\times[-c_k,c_k]\times [k,k]} \left| f(x,y,z)\right|\,,
        \end{aligned}
    \end{equation*}
    for some \(c_k\), as \(u_{k,n}\) is smooth and \([0,1]\) is closed.
    Then, we have
    \begin{equation}\label{step2:3}
        \begin{aligned}
              & \int_{0}^{1} f(x, u_{k,n}(x), u'_{k,n}(x)) \, dx \ \xrightarrow{n\ra \infty} \int_{0}^{1} f(x, u_{k}(x), u'_{k}(x)) \, dx \,,
        \end{aligned}
    \end{equation}
    so the difference becomes arbitrarily small.
    
    \item \label{step3} \textbf{Correct boundary values.}
    
    {
   Without loss of generality we can assume that \( u(0) = 0 \) and \( u(1) = 1 \). If this is not the case, we can adjust the function by adding \( u(0) + x \big( u(1) - u(0) \big) \) to the following definition of \( \phi_{n,k}\), and the calculations remain valid under this adjustment.
    }
    Let us take 
    \begin{equation*}
        w_{k,n}(x)\coloneqq u_{k,n}'(x) + \big(u_{k,n}\left(0\right)+1- u_{k,n}\left(1\right) \big) \,,
    \end{equation*}
    and define our final function
    \begin{equation}\label{step3:1}
        \phi_{k,n}(x)\coloneqq  \int _0^x w_{k,n}(t)\, dt\,.
    \end{equation}
    This allowed us to fix the boundary values \(\phi_{k,n}(0)=0\) and \(\phi_{k,n}(1)=1\). Finally, we need to check the remaining two conditions.
    
    \textbf{Condition~\ref{Lav1:2}:}\\
Since \(u_k\) differed from \(u\)  by no more than  \(\varepsilon\), we obtain the following estimate  \begin{equation}\label{step3:2}
        \begin{aligned}
            |\phi_{k,n}(x)  - u_{k,n}(x)| &=\Bigg| \int _0^x w_{k,n}\, dt            
             -\int _0^x u_{k,n}'(t)\, dt - u_{k,n}(0) \,\Bigg|\\
                       =&\,\Bigg| \int _0^x u_{k,n}(0)+1- u_{k,n}(1) \, dt- u_{k,n}(0)\Bigg|        \\
                      =&\left |(x -1) u_{k,n}(0) + x(1-u_{k,n}(1)) \right| \\
                      \leq&\left|u_{k,n}(0)\right| + |1-u_{k,n}(1)| \leq \,4\ve\,.
        \end{aligned}
    \end{equation}

    \textbf{Condition~\ref{Lav1:3}:}\\
   Firstly, we will break this into smaller parts
    \begin{equation}\label{step3:3}
        \begin{aligned}
              & \left|\int_0^1 f(x,u_{k,n}(x),u_{k,n}'(x))\, dx -\int_0^1 f(x,\phi_{k,n}(x),\phi_{k,n}'(x))\, dx \right|                            \\
              & \quad\leq \left|\int_0^1 f(x,u_{k,n}(x),u_{k,n}'(x))\, dx -\int_0^1 f(x,u_{k,n}(x),\phi_{k,n}'(x))\, dx \right|                     \\
              & \quad \ \ +\left|\int_0^1 f(x,u_{k,n}(x),\phi_{k,n}'(x))\, dx -\int_0^1 f(x,\phi_{k,n}(x),\phi_{k,n}'(x))\, dx \right|\\&\hspace{10cm} \eqqcolon I_1 + I_2 \,.
        \end{aligned}
    \end{equation}

Let us begin by analyzing \(I_1\). Using continuity and estimates similar to those in \eqref{step3:2}, we observe that:

    \(\forall_{\delta>0} \quad \exists_{k_0}\quad\forall_{k>k_0}\quad \exists_{n_k} \quad\forall_{n\geq n_k} \quad\)
    \begin{equation*}
        \left|u_{k, n}(0) + 1 - u_{k, n}(1)\right| \leq  \delta \quad \text{for all } x \in [0,1]\,.
    \end{equation*}

Next, consider the mapping \(\psi \colon x \mapsto \left( x, u_{k,n}(x), u_{k,n}'(x) \right)\).  Due to the fact that \( u_{k,n} \) is a Lipschitz continuous function of constant \(k\), it is bounded on \([0,1]\) by
\[
\psi([0,1]) \subset [0,1] \times [-C, C] \times [-k, k]\,.
\]

By the continuity of \( f \), we know that the map \( (x, y, z) \mapsto f(x, y, z) \) is uniformly continuous on \(\psi([0,1]) \subset [0,1] \times [-C, C] \times [-k, k]\). Denoting the modulus of continuity by \( \omega \), we can estimate
    \begin{multline*}
        \left| f\Big(x,u_{k,n}(x), u'_{k,n}(x)\Big) \right.\\\left.-f\Big(x,u_{k,n}(x),
        u_{k,n}'(x) + u_{k,n}(0)+1 - u_{k,n}(1)\Big)\right|\\
        \leq\omega\big|u_{k,n}(0)+1 - u_{k,n}(1) \big|\leq \omega\, \delta\,.
    \end{multline*}
    Thus, we can estimate
    \[ I_1 \leq \omega \delta\,, \]
    which is arbitrarily small. It remains to show that \(I_2\) is negligible. Indeed,
    \begin{multline*}
         I_2=  \left|\int_0^1 f(x,u_{k,n}(x),\phi_{k,n}'(x))\, dx -\int_0^1 f(x,\phi_{k,n}(x),\phi_{k,n}'(x))\, dx \right|                                          \\
            \leq  \left| \sup_{[0,1]\times \R^2}\frac{\partial f}{\partial y} (x,y,z)\big( u_{k,n} - \phi_{k,n}\big)\right|\leq M\big \|u_{k,n} - \phi_{k,n} \big\|_{\infty}\leq M\ve \,.
    \end{multline*}

    \textbf{Final verification:} \\
Let us show that the function $\phi_{k,n}$ from \eqref{step3:1} satisfies all conditions of Lemma~\ref{Lav1}. Note that it agrees with $u$ on the boundary, so we only need to check the last two conditions.

Let us verify \textbf{Condition~\ref{Lav1:2}} from Lemma~\ref{Lav1} for \(u\). From \eqref{step1:2}, \eqref{step2:2}, and \eqref{step3:2}, we obtain\\
    \(\forall_{\varepsilon>0} \quad \exists_{k_0}\quad\forall_{k>k_0}\quad \exists_{n_k} \quad\forall_{n\geq n_k} \quad \forall_{x} \quad\)
    \begin{equation*}
        \begin{aligned}
            |\phi_{k,n}(x) - u(x)| \leq\;  |\phi_{k,n}(x) - u_{k,n}(x)|&  + |u_{k,n}(x)- u_{k}(x)| \\ &+ |u_k(x) - u(x)|\leq3\ve\,.
        \end{aligned}
    \end{equation*}

    \textbf{Condition~\ref{Lav1:3}:} From \eqref{step1:3}, \eqref{step2:3}, and \eqref{step3:3}, we get\\
    \(\forall_{\varepsilon>0} \quad \exists_{k_0}\quad\forall_{k>k_0}\quad \exists_{n_k} \quad\forall_{n\geq n_k} \quad\)
    \begin{multline*}
          \left| \int_0^1 f(x,\phi_k(x),\phi_k'(x)) \,dx - \int_0^1 f(x,u(x),u'(x)) \,dx\right|                  \\
               \leq\left| \int_0^1 f(x,\phi_k(x),\phi_k'(x)) \,dx - \int_0^1 f(x,u_{k,n}(x),u_{k,n}'(x)) \,dx\right|  \\
               +\left| \int_0^1 f(x,u_{k,n}(x),u_{k,n}'(x)) \,dx - \int_0^1 f(x,u_{k}(x),u_{k}'(x)) \,dx\right|       \\
               +\left| \int_0^1 f(x,u_{k}(x),u_{k}'(x)) \,dx - \int_0^1 f(x,u(x),u'(x)) \,dx\right| \leq 3\varepsilon \,.
    \end{multline*}
    This proves the lemma.   \hfill\qedsymbol{}
\end{enumerate}

From the proof of Lemma~\ref{Lav1}, it turns out that the hypothesis might be strengthen and additionally yields convergence a.e. of the derivatives.{
\begin{remark}%
    Let $f$ and $u$ be as in Lemma~\ref{Lav1}. Then for every $\varepsilon>0$, there exists $\phi\in C^\infty([0,1])$, such that conditions~\ref{Lav1:1}--\ref{Lav1:3} are satisfied and the following holds true
		\begin{enumerate}
  \item[(L4)] $\int _0 ^1 |u'(x)-\phi'(x)|\,dx<\varepsilon\,.$ 
  \end{enumerate}
\end{remark}

\begin{proof} Fix $\varepsilon>0$ and assume that $x$ is such that $u'(x)$ is well-defined. It is enough to prove that $\{\phi_{k,n}\}$ from the proof above, see \eqref{step3:1}, satisfies
\begin{equation}\label{u'-phi'-closeness}
    \int _0^1 |u'(x)-\phi_{k,n}'(x)|\,dx<\varepsilon
\end{equation}
for all $k,n$ large enough. Let us observe that
\[|u'(x)-\phi'_{k,n}(x)|\leq |u'(x)-u'_{k}(x)|+|u_k'(x)-u'_{k,n}(x)|+|u'_{k,n}(x)-\phi'_{k,n}(x)|\,,\]
where $\{u_k\}$ is defined in~\ref{step1}, $\{u_{k,n}\}$ is defined in~\ref{step2},  whereas $k,n$ are sufficiently large as described in the proof of Lemma~\ref{Lav1}. By the very definition of $u_k$, cf. \eqref{step1:2} in~\ref{step1}, we obtain that
\[\int_0^1 |u'(x)-u'_k(x)| \,dx \leq      \int_{\{|u'(t)|>k\}} |u'(t)| \,dt \xrightarrow{k\ra \infty}0 \,.\]
{
The convergence of \(u_{k,n}\) defined in step~\ref{step2} is verified by following estimates:
}
\begin{equation*}
\begin{aligned}
        | u_k'(x) - u_{k,n}'(x)| =& \left| u_k'(x) - \frac{d}{dx} \int_{-\frac{1}{n}}^{\frac{1}{n}} u_k(x-y) \eta_{\frac{1}{n}}(y) \, dy \right| \\
    =& \left| u_k'(x) - \int_{-\frac{1}{n}}^{\frac{1}{n}} u_k'(x-y) \eta_{\frac{1}{n}}(y) \, dy \right| \\
    \leq& \left| \int_{-\frac{1}{n}}^{\frac{1}{n}} (u_k'(x) - u_k'(x-y)) \eta_{\frac{1}{n}}(y) \, dy \right| \\
    \leq& \int_{-\frac{1}{n}}^{\frac{1}{n}} \left| k y \right| \eta_{\frac{1}{n}}(y) \, dy < \frac{2k}{n}\,.
\end{aligned}
\end{equation*}
    We can change the order of differentiation and  integration by invoking the Dominated Convergence Theorem, since \( u_k \) is a Lipschitz function. In~\ref{step3}, by the definition of \(\phi_{k,n}\), we can estimate the following for \(x\in [0,1]\):
\begin{equation*}
    \left| \phi'_{k,n}(x) - u_{k,n}'(x) \right| = \left|  u_{k,n}(0) +1 - u_{k,n}(1)  \right| 
    \leq \left| u_{k,n}(0) +1 - u_{k,n}(1) \right| ,
\end{equation*}
where the right-hand side is arbitrarily small, cf. \eqref{step3:2}. Inserting the above observations to~\eqref{u'-phi'-closeness} and integrating, we conclude the proof.
\end{proof}}
\section{Consequences of the Approximation Lemma}
From Lemma~\ref{Lav1}, we can derive two important theorems:
{ 
\begin{theorem}[Lavrentiev 1927, \cite{Mascon}] \label{coro1}
    Suppose \( f \) is as in 
    Lemma~\ref{Lav1}, then the infimum (supremum) of the integral \(\int_{0}^{1} f(x, u(x), u'(x))\,dx\) in the class \( AC_*([0,1]) \) is equal to the lower (upper) bound of the same integral in the class \( C_*^{\infty}([0,1]) \).
\end{theorem}
}
\begin{proof}
Since \(  C^\infty_*([0,1]) \) is a subset of \( AC_*([0,1]) \), we immediately know that
\begin{equation*}
    \inf_{u \in AC_*([0,1])} \int_0^1 f(x,u(x),u'(x)) \, dx \leq \inf_{u \in C^\infty_*([0,1])} \int_0^1 f(x,u(x),u'(x)) \, dx \,.
\end{equation*}
It suffices to prove the reverse inequality. Let \(\{u_n\} \subset AC_*([0,1])\) be a minimizing sequence. Then for every \(\varepsilon > 0\), we can find \(n_\varepsilon\) such that for all \(n > n_\varepsilon\), it holds that
\[
\int_0^1 f(x, u_n(x), u_n'(x)) \, dx \leq \inf_{u \in AC_*([0,1])} \int_0^1 f(x, u(x), u'(x)) \, dx + \varepsilon\,.
\]
\noindent
We fix an arbitrary \(\varepsilon > 0\) and choose \(n > n_\varepsilon\). By Lemma~\ref{Lav1}, we obtain \(\varphi_n \in C^\infty_*([0,1])\) such that
\[
\left|\int_0^1 f(x, u_n(x), u_n'(x)) \, dx - \int_0^1 f(x, \varphi_n(x), \varphi_n'(x)) \, dx \right|< \varepsilon\,.
\]
\noindent
Combining the inequalities, for all \(n > n_\varepsilon\), we have
\begin{align*}
    \int_0^1 f(x,\varphi_n(x),\varphi_n'(x)) \, dx &\leq  \int_0^1 f(x, u_n(x), u_n'(x)) \, dx +\ve \\ 
        &\leq \inf_{u \in AC_*([0,1])} \int_0^1 f(x, u(x), u'(x)) \, dx + 2\varepsilon\,.
\end{align*}
As \(\varepsilon > 0\) was arbitrary, we can pass to the limit and conclude the proof.
\end{proof}

\begin{theorem}[Lavrentiev 1927, 
\cite{Mascon}]\label{coro2} 
    If there exists a map in the class \( C^{1}_* ([0,1]) \) that minimizes the integral \(\int_{0}^{1} f(x, u(x), u'(x))\,dx\) (where \(f \) satisfies the conditions of Theorem~\ref{Lav1}), then there also exists a minimum for this integral in the class \( AC_*([0,1]) \), and this minimum is attained by the same map.
\end{theorem}

\begin{proof}{
 
    Since \( C^\infty _*([0,1]) \subset  C^1_*([0,1]) \subset AC_*([0,1]) \), and Theorem~\ref{coro1} establishes that the lower bound in \(  C^\infty_*([0,1])\) and \( AC_*([0,1])\) is the same, it follows that the function above is also a minimizer of the integral in the class  \( AC_*([0,1]) \).}
\end{proof}

\appendix
\section{Supplementing Lavrentiev's proof}\label{Appendix_A}

As discussed in Section~\ref{Simple_counterexample}, here we focus on the issue of constructing a partition into a finite number of non-intersecting intervals for any absolutely continuous function, as required in point~\ref{org:step4} To address this, we provide a methodical construction of such intervals, emphasizing the critical role played by measure regularity.
{ 
We shall use the following Lemma.
\begin{lemma}\label{Lemma0}
    In a separable metric space \((X, d)\), given a nonatomic measure \(\mu\), the set of isolated points has measure zero.
\end{lemma}
As a consequence of Lemma~\ref{Lemma0}, we have the following fact.
\begin{lemma}\label{remark:1}
	In a separable metric space \((X,d)\), given a nonatomic measure \(\mu\), consider a measurable set \(A\). If \(\mu\) is inner regular, the following holds:
	\[
		\mu(A) = \sup \big\{ \mu(F) \colon F \subseteq A, F \text{ is  {perfect}} \big\}\,.
	\]
\end{lemma}
}
{  The following lemma provides the foundation for our construction.
\begin{lemma}\label{lem:perfect_set_partition}
    Let \( P_1 \) and \( P_2 \) be two disjoint perfect subsets of the interval \([0,1]\). Then, there exists a finite collection of disjoint intervals \( \mathcal{B} = \{B_1, B_2, \dots, B_n\} \) such that:  
    \begin{enumerate}  
        \item \( \bigcup_{i=1}^{n} B_i = [0,1] \),  
        \item For each interval \( B_i \), either  
              \(
              P_1 \cap B_i = \emptyset \) or \( P_2 \cap B_i = \emptyset,
              \)
        \item\label{itm:3:perf} Each interval is of the form \([a,b)\) for \( B_i \), \( i \in \{1,2,\dots,n-1\} \), and \( B_n \) is closed.
    \end{enumerate}  
\end{lemma}  

\begin{proof}  
Define
\[
y_x \coloneqq
\begin{cases}
\inf \{ y > x \colon y \in P_2 \} & \text{if } x \in P_1\,, \\
\inf \{ y > x \colon y \in P_1 \} & \text{if } x \in P_2\,,
\end{cases}
\] 
and, if \( 0 \notin P_1 \cup P_2 \), set
\[
y_0 \coloneqq \inf \{ y > 0 \colon y \in P_1 \cup P_2 \} \,.
\]
For each \( x \in (P_1 \cup P_2) \cup \{0\} \) set
\[
I_x \coloneqq 
\begin{cases}
[x, y_x) & \text{if } y_x<\infty\,, \\
[x, 1] & \text{otherwise.}
\end{cases}
\]
This construction yields a family of interval \( \{ I_x \} \), each intersecting at most one of the sets \( P_1 \) or \( P_2 \).
Next, remove from this collection any interval strictly contained in another and denote the resulting disjoint collection by \( \mathcal{I} \). These intervals are not degenerate, since \( P_1 \) and \( P_2 \) are disjoint closed sets.

We now prove that \( \mathcal{I} \) is finite. Suppose, for contradiction, that there are infinitely many intervals in \( \mathcal{I} \) with left endpoints in \( P_1 \) (the same argument works symmetrically for \( P_2 \)). Then we can choose a strictly monotone sequence \( \{a_i\}_{i=1}^\infty \subset P_1 \) of such left endpoints converging to some \( a \in P_1 \), by the closedness of \( P_1 \).
By construction, between \( a_i \) and \( a_{i+1} \) there must exist a point \( b_i \in P_2 \).
For instance, if the sequence \( \{a_i\} \) is increasing, one may take \( b_i = y_{a_i} \); if it is decreasing, then \( b_i = y_{a_{i+1}} \).
Hence, we obtain a sequence \( \{b_i\}_{i=1}^\infty \subset P_2 \) with \( b_i \in (a_i, a_{i+1}) \) (or \( b_i \in (a_{i+1}, a_{i}) \) if the sequence is decreasing), and thus \( \lim_{i \to \infty} b_i = a \in P_1 \). This contradicts the assumption that \( P_1 \cap P_2 = \emptyset \).

Therefore, the collection \( \mathcal{I} \) must be finite, completing the proof.
\end{proof}
}
To construct the required intervals, instead of searching for the sets mentioned in point~\ref{org:step3} of the sketch of the original proof, we can begin directly by identifying the family described in point~\ref{org:step4}
\begin{lemma}

    For every \(\ve > 0\), there exists a family \(\mathcal{B} = \{B_1, B_2, \ldots, B_n\}\) of disjoint intervals and a set \(P_\ve\subseteq[0,1]\) such that:
    \begin{enumerate}
        \item \label{cond1} \(\lambda(P_\ve) > 1 - \ve\) and \(P_\ve\) is perfect,
        \item \label{cond2} \(\bigcup B_i = [0, 1]\,,\)
        \item \label{cond3}for each \(i\) we have \(\diam\left(u'\left(B_i \cap P_\ve\right)\right) \leq \ve\), where \(\diam(\cdot)\) denotes the diameter of a set.
    \end{enumerate}
\end{lemma}
\begin{proof} 
Then, the sets \(A_i\) are simply the intersections of \(P_\ve\) and \(B_i\). Let us prove that such a family exists.
	Choose \(\ve>0\). Let \(n\) be so large that \(\frac{1}{n}<\ve\). For each \(k\in \mathbb{Z}\), define
 \[
 C_k\coloneqq\left(u'\right)^{-1}\Bigg(\left[\frac{k}{n}, \frac{k+1}{n}\right)\Bigg)\,.
 \]
 Those sets are measurable as the preimage of a measurable function and pairwise disjoint. Since \(u'\in L^1\left([0,1]\right)\), there exists \(M\in \N\) such that:
	\[\exists_{M\in \mathbb{N}} \quad \lambda \left( \bigcup _{k<-M}C_k \cup \bigcup _{k>M} C_k\right)\leq \frac{\ve}{4}\,.\]
	We fix such $M$ and we notice that it holds
	\[\lambda \left( \bigcup _{k=-M}^{M}C_k\right)\geq 1- \frac{\ve}{4}\,.\]
	
	\noindent
	From the inner regularity of the Lebesgue measure we can find compact sets \( P_{\ve,k}\) such that:
	\[\forall_{k\in [-M,M]\cap \mathbb{Z}} \quad \exists_{P_{\ve,k} \in C_k} \quad \lambda(C_k\setminus  P_{\ve,k}) \leq \frac{\ve}{4M}\,.\]
	Without loss of generality, we can assume that \(P_{\ve,k}\) are perfect sets as per Lemma~\ref{remark:1}. Let \[  P_\ve^m \coloneqq \bigcup
 _{k=-M}^M   P_{\ve,k}\,.\]Then, we have
\begin{equation}\label{sup:con1}
\lambda ( P_\ve) >1- \frac{(2M+1)\ve}{4M} - \frac{\ve}{4}\,.
	\end{equation}
    { 
  Since \( P_{\varepsilon} \) is a finite union of closed sets, it is also closed. Define the sets
    \[
    F_m \coloneqq \bigcup_{k=-M}^{m} P_{\varepsilon,k}, \quad G_m \coloneqq \bigcup_{k=m+1}^{M} P_{\varepsilon,k}, \quad \text{for } m \in [-M, M-1] \cap \mathbb{Z}\,.
    \]
    Using Lemma~\ref{lem:perfect_set_partition} on the pairs \( (G_m, F_m) \), we obtain finite families of intervals \(\mathcal{B}'_m\). Let \(\mathcal{B'} = \bigcup_{m=-M}^{M-1} B'_m \) be sum of those families.
 }
To construct our family \(\mathcal{B}\), take the endpoints of every interval in \(\mathcal{B'}\) and sort them to get the sequence \(\{a_n\}_{n=1}^N\). Set \(a_0 = 0\), then let
\[
\mathcal{B} = \big\{[a_i, a_{i+1}) \colon i \in \left\{0, 1, \ldots, N\right\}\big\} \cup \left\{\left[a_{N}, 1\right]\right\}.
\]
In~\eqref{sup:con1}, we ensured that condition~\ref{cond1} is satisfied. This construction also guarantees that the sum of these intervals covers the entire interval \([0,1]\), thus satisfying condition~\ref{cond2}. 

{
Now, let us check condition~\ref{cond3}. Consider any interval \(B\in\mathcal{B}\). If this interval has an empty intersection with \(P_\varepsilon\), then there is nothing to prove, as \(\diam(\Oset) = 0 < \varepsilon\). 
If the interval has a non-empty intersection with \(P_\varepsilon\), it must intersect \(P_{\varepsilon,k}\) for some \(k\). 
Suppose it also intersects \(P_{\varepsilon,j}\) for some \(j > k\). Then, \(\mathcal{B}'_k\) is a partition of \([0,1]\) such that every interval from this family intersects exactly one of the sets \(P_{\varepsilon,k}\) or \(P_{\varepsilon,j}\). Since the interval \(B\) is contained in some interval from \(\mathcal{B}'_k\), this yields a contradiction.
}
Therefore, the diameter \(\diam\left(u'\left(B_i \cap P_\varepsilon\right)\right) = \diam\left(u'\left(B_i \cap P_{\varepsilon, k}\right)\right) \leq \varepsilon\), as \(P_{\varepsilon, k} \subset C_k\).
\end{proof}

\medskip
\textbf{Acknowledgements:}
The author would like to thank Michał Borowski  and Prof.~Witold Marciszewski for the insightful discussions on the new proofs. 
Appreciation is also due to friends Łukasz Wodnicki, Tomasz Jurczak, and Natalia Olszewska for their support and for reading drafts of the Lavrentiev example proof.
Special thanks to dr hab. Iwona Chlebicka for her invaluable guidance and support throughout this work.

\bibliographystyle{abbrv}

\end{document}